\documentclass{article}
\usepackage{amsfonts}
\usepackage{amssymb}
\usepackage{amsmath}
\usepackage[colorlinks=true]{hyperref}
\usepackage{color}

\begin{document}

\title{Two hierarchies of spline interpolations. Practical algorithms for multivariate higher order splines.}
\author{Cristian Constantin Lalescu\footnote{Statistical and Plasma Physics, Universit\'e Libre de Bruxelles, Campus Plaine, CP 231, B-1050 Brussels, Belgium} \footnote{\textit{Electronic address:} \texttt{clalescu@ulb.ac.be}}}

\maketitle

\begin{abstract}
A systematic construction of higher order splines using two hierarchies of polynomials is presented.
Explicit instructions on how to implement one of these hierarchies are given.
The results are limited to interpolations on regular, rectangular grids, but an approach to other types of grids is also discussed.
\end{abstract}
\vfill
\tableofcontents
\newpage

\section{Introduction}

The purpose of this work is to construct smooth interpolants for functions that are only known on the nodes of a regular rectangular grid.
Cubic splines have recently been used in the context of particle or virtual particle tracking in complex fields, \cite{Mackay, Lekien_tricubic, Homann_impact}.
A more elaborate discussion on integrating particle trajectories in interpolated fields, and the advantages of using higher order splines, is to be published elsewere\footnote{C.C. Lalescu, B. Teaca, and D. Carati, ``Implementation of high order spline interpolations for tracking test particles in discretized fields'', submitted for publication.}.
In this work just a review of the construction of the spline interpolations is presented, without rigurous proofs.

Consider the functions of $D$ variables
\begin{equation}
f:\mathbb{R}^D \rightarrow \mathbb{R}
\end{equation}
and assume they can only be computed on the grid
\begin{equation}
G = \prod_{j=1}^D h_j \mathbb{Z},
\end{equation}
where $0<h_j<1$ are the grid constants and $h\mathbb{Z} = \{h z \big| z \in \mathbb{Z}\}$; a grid cell is defined as:
\begin{equation}
C_{(z_1,z_2,\dots,z_D)} \equiv \prod_{j=1}^D [h_jz_j, h_j(z_j+1)]
\end{equation}

Basically, the values $f(x)$ are only available if $x\in G$ --- this can happen for a variety of reasons, the simplest being that they are measured from an experiment.
In the following we will say that we are computing approximations of $f$.
In a more rigurous setting we would say that we are building an array of spline functions that converge to a certain limit under certain conditions; when the function $f$ is sufficiently nice (a term that still needs a clear definition), it will be equal to that limit.
For instance, one of the properties of this limit is that any of its Taylor expansions converge everywhere (for example any function with a finite discrete Fourier representation).

A polynomial spline is a function defined on $\mathbb{R}^D$, that is a polynomial on each cell and it has continuous derivatives \emph{everywhere} up to a certain order.
Note that the set of grid constants $\{h_j\}$ is a given, and the limit spoken of before is the limit of very large orders of the polynomials entering the spline.
This limit will always exist as long as $f$ is well defined on the grid $G$; it is true that for certain functions the limits for different grids will not be equal, but these functions are irrelevant here.

\section{One dimensional case} \label{sec: splines}

\subsection{Hermite splines}

Assume that the correct values of a function $f(x)$ are known on $h\mathbb{Z}$, and we are interested in finding an approximation for the interval $[x_0,x_0 + h]$ (the cell $C_{(x_0)}$).
Without loss of generality, we express the variable $x$ in $h$ units, and the formulas will be deduced for the interval $[0,1]$ (the cell $C_{(0)}$).

On $[0,1]$ we construct the $n$-th order polynomial
\begin{equation}
s^{(n)}(x) = \sum_{k=0}^n a^{(n)}_k x^k
\end{equation}
For Hermite spline interpolation, it is imposed that, on the enclosing grid nodes $s^{(n)}$ coincides with the original function and the derivatives of $s^{(n)}$ up the order $m\equiv (n-1)/2$ coincide with the derivatives of the original function:
\begin{equation}
\left[\tfrac{d^l s^{(n)}}{dx^l}(x) = f^{(l)}(x)\right]_{ x\in\{0,1\}},\ l=\overline{0,m}
\label{eq:spline_eq}
\end{equation}
where we called the $l$-th order derivative $f^{(l)}\equiv \tfrac{d^l}{dx^l} f$.

By solving the linear system of equations \eqref{eq:spline_eq} the coefficients $a^{(n)}_k$ depending on $f(0),f(1),f'(0)\dots$ can be easily found:
\begin{equation}
a^{(n)}_k = \sum_{l=0}^m \sum_{i=0}^1 b^{(n)}_{kli} f^{(l)}(i) \label{ak} \, .
\end{equation}
Here we will discuss a method that avoids the computation of these coefficients.
This leads to rewriting the expression of the spline as:
\begin{align}
s^{(n)}(x) &= \sum_{k=0}^n \left(\sum_{l=0}^m \sum_{i=0}^1 b^{(n)}_{kli} f^{(l)}(i)\right) x^k \crcr
&= \sum_{l=0}^m \sum_{i=0}^1 f^{(l)}(i) \left( \sum_{k=0}^n b^{(n)}_{kli} x^k \right) \crcr
s^{(n)}(x) &= \sum_{l=0}^m \sum_{i=0,1} f^{(l)}(i) \alpha^{(n,l)}_i (x)
\label{eq:1D_Hermite}
\end{align}
where the $\alpha$ polynomials can be found for specific values of $n$ by symbolic computation; they can be defined equivalently as the solutions of the following system of equations:
\begin{equation}
\tfrac{d^l}{dx^l}\alpha^{(n,l_0)}_i (x) = \delta_{i,j} \delta_{l_0,l}\Big|_{ j,i\in\{0,1\},x=j},\ l=\overline{0,m}
\label{eq:alpha_i_equations}
\end{equation}

The set $\{s^{(n)}\}$ is a hierarchy of spline approximations (Hermite splines), and it is based on the \emph{spline polynomials of the first kind} $\alpha^{(n,l)}_i$.
These polynomials have the following explicit expressions:
\begin{align}
\alpha^{(n,l)}_0 (x) &= \frac{x^l}{l!}(1-x)^{m+1} \sum_{k=0}^{m-l} \frac{(m+k)!}{m! k!}x^k \\
\alpha^{(n,l)}_1 (x) &= \frac{(x-1)^l}{l!}x^{m+1} \sum_{k=0}^{m-l} \frac{(m+k)!}{m! k!}(1-x)^k
\end{align}
(with $\alpha^{(n,l)}_1 (x) = (-1)^l \alpha^{(n,l)}_0 (1-x)$); see appendix \ref{sec:spline_polynomials} for a proof.

The distance $||s^{(n+1)} - s^{(n)}||$ will obvioulsy go to $0$ for large $n$, because the splines coincide with the Taylor approximation around the two nodes up to the order $m$; note that this distance can be defined as the maximum of the distance $|s^{(n+1)}(x) - s^{(n)}(x)|$.
And there is a class of functions $f$ that are limits of such approximations (all piecewise polynomials for instance).

\subsection{Grid splines}

If centered differences are used to compute the derivatives, the formula can be further adapted.
The result is called here a \emph{grid spline} (it should be a type of b-spline), as it is found solely from the values of the function on the grid.
Centered differences are used because generally for $D>1$ dimensions --- when discussing practical implementations of higher order splines --- the memory cost of keeping all the derivatives necessary is prohibitive (they are as many as the coefficients $a^{(n)}_k$).
A finite difference (for this rescaled variable) is just a linear combination of the numbers $f(z)$ with $z\in\mathbb{Z}$.
Also, it is crucial that the use of \emph{centered} differences imposes the continuity of the derivatives when the polynomials are put together.

To be more specific, rewrite the Hermite spline as
\begin{equation}
s^{(n)}(x) = \sum_{l=0}^m \sum_{i=0,1} t_i^{(l)} \alpha^{(n,l)}_i (x).
\end{equation}
The smoothness of the interpolant is only determined by the fact that the coefficients of the Taylor expansions used, $t_i^{(l)}$, are determined by the grid node $i$ (they are the same whether the node is approached from the left or from the right, or more simply they are invariant to the cell).
In practice the centered differences are the coefficients of the Lagrange interpolation polynomial $\sum_{k=0}^{2g} t_{(i,g)}^{(k)} x^k$, with the coefficients determined from the equations
\begin{equation}
\left[\sum_{k=0}^{2g} t_{(i,g)}^{(k)} x^k = f(i+x)\right]_{x=\overline{-g,g}}.
\end{equation}

Because the information contained in this Taylor expansion must be invariant to the cell, it must be obtained from a set of grid nodes that is symmetrical to the current grid node $i$ (thus the equations are solved at ${x=\overline{-g,g}}$, so a unique solution is obtained for a polynomial of order $2g$).

It is then obvious that a centered difference $f^{<q,l>}\approx f^{(l)}$ can be written as:
\begin{equation}
f^{<q,l>}(i) \equiv t_{(i,g)}^{(l)} = \sum_{k=-g}^{g} c^{(q,l)}_k f(i+k)
\end{equation}
where the coefficients $c^{(q,l)}_k$ are numbers that only depend on $g$.
This means that the final formula of the spline will use the values of the function on the $q=2g+2$ nodes $-g,\dots,0,1,\dots,g+1$:
\begin{equation}
s^{(n,q)}(x) = \sum_{i=0-g}^{1+g} f(i) \beta^{(n,q)}_i (x)
\label{eq:general_spline}
\end{equation}
where the \emph{spline polynomials of the second kind} $\beta^{(n,q)}_i(x)$ can be found for a fixed $n$ and $q$ by symbolic computation.

We give the $\beta^{(5,4)}_i$-s as an example:
\begin{align}
\beta^{(5,4)}_{-1}(x) &= \frac{1}{2}{\left( x-1\right) }^{3}\,x\,\left( 2\,x+1\right)\\
\beta^{(5,4)}_0(x) &= -\frac{1}{2}\left( x-1\right) \,\left( 6\,{x}^{4}-9\,{x}^{3}+2\,x+2\right)\\
\beta^{(5,4)}_1(x) &= \frac{1}{2}x\,\left( 6\,{x}^{4}-15\,{x}^{3}+9\,{x}^{2}+x+1\right)\\
\beta^{(5,4)}_2(x) &= -\frac{1}{2}\left( x-1\right) \,{x}^{3}\,\left( 2\,x-3\right)
\end{align}

In terms of accuracy, using the distances $||f^{<q,l>} - f^{(l)}||$ one should be able to find the distance $||s^{(n,q)} - s^{(n)}||$.

As a sidenote, it is quite clear that the same approach can be used for irregular grids (the Taylor expansions can still be approximated using adjacent nodes).
However, some of the simplifications that can be made for regular grids will no longer be possible.

\section{D dimensions}

Once the spline polynomials of the first and the second kind are available, the spline interpolation formula can be directly generalized to a scalar function of $D$ variables (expressed in $h_j$ units, and shifted to $[0,1]^D$):
\begin{equation}
s^{(n)}(x_1,x_2,\dots,x_D) = \sum_{l_1,\dots,l_D=0}^m \sum_{i_1,\dots,i_D=0,1} f^{(l_1,\dots,l_D)}(i_1,\dots,i_D) \prod_{k=1}^D \alpha^{(n,l_k)}_{i_k} (x_k)
\label{eq:Hermite_spline_DD}
\end{equation}
\begin{equation}
s^{(n,q)}(x_1,x_2,\dots,x_D) = \sum_{i_1,\dots,i_D=0-g}^{1+g} \Big(f(i_1, \dots,i_D ) \prod_{j=1}^D \beta^{(n,q)}_{i_j}(x_j)\Big)
\label{eq:grid_spline_DD}
\end{equation}
A rather interesting observation is that for the Hermite splines exactly $2^D (m+1)^D$ input values $f^{(l_1,\dots,l_D)}(i_1,\dots,i_D)$ are required for a given $n=2m+1$, while for grid splines $q^D$ input values are required for a given $q$ (and $m \leq 2g$, thus $n\leq 2q-3$).
This means that grid splines generally achieve a given degree of smoothness from less information than a Hermite spline --- and the price is probably a much larger error.
Also, the number of input values is the number of terms in the sum to be computed, thus grid splines will be faster to compute than the corresponding Hermite splines (up to a maximum $q$ that will depend on the order).
As a final note, the mixed derivatives of these fields are also smooth:
\begin{equation}
\left(\prod_{j=1}^N \left(\frac{\partial}{\partial x_j}\right)^{l_j}\right) s^{(n)}, \left(\prod_{j=1}^N \left(\frac{\partial}{\partial x_j}\right)^{l_j}\right) s^{(n,q)} \in \mathcal{C}^{m-\max\{l_j\}}
\end{equation}
(whenever $m \geq \max\{l_j\}$).

\section{Implementation of the grid splines}

In practice the Hermite splines will probably not be very useful, as they require too much memory.
Other than that, their implementation should be similar to that of the grid splines.
Note that we mention ``parallelized codes'' in the following; this refers specifically to cases of computer programs that run on several processors at once, with the memory divided between them, and these programs work with physical fields, each processor keeping a slice of these fields in its memory.

\subsection{Full expressions of the 1D grid splines}

You will need to use a computer algebra system.
For $q=4$, define the following functions:
\begin{align}
t_{(0,1)}(h) &= t^0_{(0,1)} + t^1_{(0,1)} h + t^2_{(0,1)} \frac{h^2}{2!} \\
t_{(1,1)}(h) &= t^0_{(1,1)} + t^1_{(1,1)} h + t^1_{(1,1)} \frac{h^2}{2!} \\
s^{(3,4)}(x) &= s^{(3,4)}_0 + s^{(3,4)}_1 x + s^{(3,4)}_2 x^2 + s^{(3,4)}_3 x^3 \\
s^{(5,4)}(x) &= s^{(5,4)}_0 + s^{(5,4)}_1 x + s^{(5,4)}_2 x^2 + s^{(5,4)}_3 x^3 + s^{(5,4)}_4 x^4 + s^{(5,4)}_5 x^5
\end{align}
Afterwards, solve the following systems of equations (symbolically):

\begin{equation}
\left\{
\begin{aligned}
&t_{(0,1)}(0) = f(0), &t_{(0,1)}(-1) = &f(-1), &t_{(0,1)}(1) = &f(1), \\
&t_{(1,1)}(0) = f(1), &t_{(1,1)}(-1) = &f(0), &t_{(1,1)}(1) = &f(2),
\end{aligned}
\right.
\label{eq:centered_differences_q4}
\end{equation}

\begin{equation}
\left\{\begin{aligned}
&&s^{(3,4)}(0) &= &t^0_{(0,1)}, &&s^{(3,4)}(1) &= &t^0_{(1,1)}, \\
&\biggr[\frac{d}{dx}&s^{(3,4)}(x)\biggr]_{x=0} &= &t^1_{(0,1)}, &\biggr[\frac{d}{dx}&s^{(3,4)}(x)\biggr]_{x=1} &= &t^1_{(1,1)},
\end{aligned}\right.
\end{equation}

\begin{equation}
\left\{\begin{aligned}
&&s^{(5,4)}(0) &= &t^0_{(0,1)}, &&s^{(5,4)}(1) &= &t^0_{(1,1)}, \\
&\biggr[\frac{d}{dx}&s^{(5,4)}(x)\biggr]_{x=0} &= &t^1_{(0,1)}, &\biggr[\frac{d}{dx}&s^{(5,4)}(x)\biggr]_{x=1} &= &t^1_{(1,1)}, \\
&\biggr[\frac{d^2}{dx^2}&s^{(5,4)}(x)\biggr]_{x=0} &= &t^2_{(0,1)}, &\biggr[\frac{d^2}{dx^2}&s^{(5,4)}(x)\biggr]_{x=1} &= &t^2_{(1,1)},
\end{aligned}
\right.
\end{equation}

You can either solve them one at a time (and afterwards replace the solution of \eqref{eq:centered_differences_q4} into the expressions of the splines), or you can solve them all at once.
For larger numbers of grid points $q$, the process is identical.

\subsection{Spline polynomials}

Note that the same systems of equations need to be solved for the Hermite splines, just that the centered differences $t^q_{(i,l)}$ should be replaced with the $f^{(l)}(i)$, and the spline polynomials of the first kind can be found as:
\begin{equation}
\alpha^{(n,l)}_i (x) \equiv \frac{d}{d(f^{(l)}(i))} s^{(n)}(x)
\end{equation}
For the spline polynomials of the second kind, the values of the function come into play (and the following expression makes sense after the centered differences have been introduced into the expressions of the splines):
\begin{equation}
\beta^{(n,q)}_i (x) \equiv \frac{d}{d(f(i))} s^{(n,q)}(x)
\end{equation}
These two equations are the simplest way to find the polynomials.
Considering that the full expressions of the splines are linear in the $t^l_{(i,g)}$-s which are linear in the $f(i)$-s, they are equivalent to rearanging the terms in the sum and extracting the ``coefficients'' of the $f(i)$-s, as in the definitions.

\subsection{Implementation}

After finding the spline polynomials of the second kind, they should be put into their Horner forms (for fast computation).
The following steps depend on the programmer's taste and abilities mostly, but we recommend the implementation of the subroutines that follow.
They are easy to adapt for the case of a parallelized code, and this implementation proved to be quite efficient and easy to work with (easy to expand for more cases, explain to other users, check for errors when necessary).
Note that the same structure can be used for interpolating derivatives of the field if necessary, just that subroutines for computing the derivatives of the spline polynomials have to be added.

Recommended structure of the interpolation code (subroutines):
\begin{enumerate}
\item \underline{spline polynomials $(n,q)$}
\begin{itemize}
\item input: fraction $\xi$
\item output: $\gamma_i = \beta^{(n,q)}_i (\xi)$ (array of dimension $q$)
\end{itemize}

\item \underline{grid coordinates}
\begin{itemize}
\item input: ``normal'' point coordinates $(x_1,\dots,x_D)$
\item algorithm: compute each $\hat x_j \equiv \lfloor x_j/h_j \rfloor$ and each $\tilde x_1 \equiv x_j/h_j - \hat x_j$.
\item output:
\begin{itemize}
\item the set of integers $(\hat x_1, \dots, \hat x_D)$
\item the set of fractions $(\tilde x_1, \dots, \tilde x_D)$
\end{itemize}
\end{itemize}

\item \underline{spline formula}
\begin{itemize}
\item input:
\begin{itemize}
\item the set of fractions $(\tilde x_1, \dots, \tilde x_D)$
\item the type of spline $(n,q)$
\item a pointer (or similar notion) $\tilde f$ to an array containing the information about the local field (the values of the field on the nodes of the cell $C_{(\hat x_1, \dots, \hat x_D)}$ and the necessary neighbouring cells), shifted such that $\tilde f(0,0,\dots,0) = f(x_1,\dots,x_D)$.
\end{itemize}
\item algorithm: compute the polynomials $\beta^{(n,q)}_i(\tilde x_j)$ in the array $\gamma_{ij}$ (by calling the spline polynomials subroutine), then compute the sum
\begin{equation}
\hat f = \sum_{i_1,\dots,i_D=0-g}^{1+g} \Big(\tilde f(i_1, \dots,i_D ) \prod_{j=1}^D \gamma_{ij}\Big);
\label{eq:grid_spline_implementedformula}
\end{equation}
note that testing shows it is more efficient to introduce as little \texttt{do while} loops as possible --- for our 3D implementation, for $q=4$ the sum is written in full in the source code, as introducing \texttt{do while} type loops slows down the code considerably.
For higher values of $q$ we just have one \texttt{do while} loop for one of the variables.
\item output: the approximation $\hat f$.
\end{itemize}

\item \underline{wrapper}
\begin{itemize}
\item input:
\begin{itemize}
\item ``normal'' point coordinates $(x_1,\dots,x_D)$
\item a pointer $f$ to the array containing the information about the entire field
\item the type of spline $(n,q)$
\end{itemize}
\item algorithm: put the local field in the array $\tilde f$ from $f$ and then compute the approximation $\hat f$ (using the above subroutines)
\item output: the approximation $\hat f$
\end{itemize}

\end{enumerate}

The wrapper is very useful in the case of a parallelized code.
All the operations related to bringing together information spread on possibly several processors can be placed inside the wrapper, allowing for easy debugging and maintenance of the code.

As an example of a parallelized version, in our implementation a 3D field is divided along the $z$ coordinate between processors.
The field is periodic in all directions, and obtaining $\tilde f$ implies a little care in regards to the $z$ coordinate, but it basically just requires a normal use of the \texttt{MODULO} operator.
We impose that each processor has at least $q$ nodes on the $z$ direction in its memory, so that the formula contains at most information from two processors (``low'' and ``up'').
We then compute \eqref{eq:grid_spline_implementedformula} on each processor, only for the nodes that are in ``its domain'', and we then sum the two resulting values $\hat f = \hat f_{low} + \hat f_{up}$; the amount of information passed between processors is thus kept to a minimum.

\appendix
\section{\appendixname: Uniqueness}

Construction for general case:
\begin{align}
t{(\mathbf{i},g)}(\mathbf{u}) &= \sum_{k_1,\dots,k_D=0}^{2g} t^\mathbf{k}_{(\mathbf{i},g)} \prod_{j=1}^D u_j^{k_j} \\
s^{(n,q)}(\mathbf{x}) &= \sum_{k_1,\dots,k_D=0}^{n} s^{(n,q)}_\mathbf{k} \prod_{j=1}^D x_j^{k_j}
\end{align}
with $i_1,\dots,i_D \in \{0,1\}$.

The Taylor expansions (the centered differences) are given by the system of equations
\begin{equation}
t_{(\mathbf{i},g)} (\mathbf{v}) = f(\mathbf{i+v})
\end{equation}
where $v_1,\dots,v_D = \overline{-g,g}$.
This system of equations has a unique solution $t^\mathbf{k}_{(\mathbf{i},g)} \approx f^{(\mathbf{k})}$ (there are as many $t^\mathbf{k}_{(\mathbf{i},g)}$ as there are $f(\mathbf{i+v})$, and this is in fact a simple Lagrange interpolation).

The splines are given by the system of equations
\begin{equation}
\left[\left(\prod_{j=1}^D \left(\frac{d}{d x_j}\right)^{l_j} \right) s^{(n,q)}(\mathbf{x}) \right]_{\mathbf{x} = \mathbf{i}} = t^\mathbf{l}_{(\mathbf{i},g)}
\end{equation}
which also has a unique solution (and $q\equiv2g+2$),
so formula \eqref{eq:grid_spline_DD} must give this unique solution.

A generalization of this method would be to go to other types of grids.
For instance one could imagine the case where in three dimensions there is a grid made up of equilateral triangles in the $(x,y)$ plane, and squares in the other two planes.
What changes is the fact that the sums are a bit more complicated.
The crucial property that has to be preserved to have a smooth approximation is that the Taylor expansion must be the same no matter from which cell we approach a given node.

\section{\appendixname: Spline polynomials}
\label{sec:spline_polynomials}
For a rigurous analysis of the errors of these approximations, the distances $||f^{<q,l>} - f^{(l)}||$, $||s^{(n,q)} - s^{(n)}||$ and $||s^{(n)} - f||$ should be computed.
After using the method presented above to find the $\alpha$ polynomials up to $n=19$, it can be seen that all these spline polynomials of the first kind have a simple form:
\begin{align}
\alpha^{(n,l)}_0 (x) &= \frac{x^l}{l!}(1-x)^{m+1} \sum_{k=0}^{m-l} \frac{(m+k)!}{m! k!}x^k \\
\alpha^{(n,l)}_1 (x) &= \frac{(x-1)^l}{l!}x^{m+1} \sum_{k=0}^{m-l} \frac{(m+k)!}{m! k!}(1-x)^k
\end{align}
(with $\alpha^{(n,l)}_1 (x) = (-1)^l \alpha^{(n,l)}_0 (1-x)$).

It would be useful to have an explicit formula for the $\alpha$ polynomials in the general case, as it would allow for a more straightforward treatment of the error $||s^{(n+1)} - s^{(n)}||$, which is close to $||s^{(n)} - f||$.

Obviously, if this form is proven to apply for $\alpha^{(n,l)}_0$, the form for $\alpha^{(n,l)}_1$ must also be correct.

First, notice that for $x$ very close to $1$ ($x=1-y$, with $y$ small), the Taylor expansion of $\alpha^{(n,l)}_0$ begins at $y^{m+1}$, and this is half of the proof.
To continue, note rewrite the polynomial as
\begin{equation}
\alpha^{(n,l)}_0 (x) = \frac{1}{m!l!} \sum_{k=0}^{m-l} \sum_{j=0}^{m+1} \frac{(m+k)!}{k!} \frac{(m+1)! (-1)^j}{j! (m+1-j)!} x^{l+k+j}
\end{equation}
For this proof the coefficient of $x^{l_0}$ in $\alpha^{(n,l)}_0$, for $0\leq l_0 \leq m$ is needed.
Obviously, for $l_0 < l$ this coefficient is 0:
\begin{align}
l+k+j &= l_0 \\
k+j &= l_0 - l\\
\text{but } k+j &\geq 0 \\
\text{so } l_0-l &\geq 0
\end{align}
For $l_0=l$, the coefficient is given by the term with $k+j=0$:
\begin{equation}
\frac{1}{l!} \frac{(m+0)!}{m! 0!} \frac{(m+1)! (-1)^0}{0! (m+1-0)!} = \frac{1}{l!}
\end{equation}
which is what is needed.
For $l_0-l=\Delta l \geq 1$ we have:
\begin{equation}
g(m,l,\Delta l) = \frac{1}{l!} \sum_{k=0}^{m-l} \sum_{j=0}^{m+1} \frac{(m+k)!}{m! k!} \frac{(m+1)! (-1)^j}{j! (m+1-j)!} \delta_{k+j,\Delta l}
\end{equation}
(with the Kronecker $\delta$).
This formula simplifies to
\begin{equation}
g(m,l,\Delta l) = \frac{1}{l!} \sum_{k=0}^{\Delta l} \frac{(m+k)!}{m! k!} \frac{(m+1)! (-1)^{\Delta l - k}}{(\Delta l - k)! (m+1-\Delta l + k)!}
\end{equation}
And in fact this sum is, from \cite{mathematica}:
\begin{equation}
g(m,l,\Delta l) = \frac{(-1)^{\Delta l} \sin (\pi \Delta l)}{l! \pi \Delta l}
\end{equation}
which is 0 for integer values of $\Delta l$ (note that $\Delta l \leq m$ for the sum to make sense).

\end{document}